\documentclass[10pt,a4paper]{article}
\usepackage[margin=1.3in]{geometry}
\usepackage[T1]{fontenc}
\usepackage{tikz}
\usetikzlibrary{positioning}
\usepackage[utf8]{inputenc}
\usepackage{graphicx}
\usepackage{mathtools}
\usepackage{amssymb}
\usepackage{amsthm}
\usepackage{thmtools}
\usepackage{xcolor}
\usepackage{nameref}
\usepackage{hyperref}
\usepackage{amsmath}
\usepackage{amsfonts}
\usepackage{wrapfig}
\usepackage{ulem}
\usepackage[autostyle]{csquotes}
\usepackage{layout}
\usepackage{tabularx}
\usepackage{authblk}
\usetikzlibrary{shapes.geometric, positioning}
\usepackage{caption}
\usepackage{subcaption}
\usetikzlibrary{positioning,chains,fit,shapes,calc}
\usepackage{tikz}
\theoremstyle{definition}
\newtheorem{defn}{Definition}[section]
\newtheorem{theorem}{Theorem}[section]

\newtheorem{lemma}[theorem]{Lemma}
\newtheorem{prop}{Proposition}[section]
\newtheorem{note}{Note}

\newtheorem{exmp}{Example}[section]
\theoremstyle{remark}
\newtheorem{remark}{Remark}[section]

\date{}
\begin{document}

\title{\sc Co-maximal Hypergraph on $D_n$ }

\author{{ \sc Sachin Ballal}\footnote{Corresponding Author}~~  and \sc{ Ardra A N} }
\affil{School of Mathematics and Statistics, University of Hyderabad,
	500046, India\\ sachinballal@uohyd.ac.in, 23mmpp02@uohyd.ac.in}

\date{}
\maketitle

\begin{abstract}
	Let $G$ be a group and $S$ be the set of all non-trivial proper subgroups of $G$. \textit{The co-maximal hypergraph of $G$}, denoted by $Co_\mathcal{H}(G)$, is a hypergraph whose vertex set  is  $\{H \in S \,\, | \,\, H K = G \,\, \text{for some} \,  K \in S  \}$ and  hyperedges are the maximal subsets of the vertex set with the property that the product of  any two vertices is equal to $G$. The aim of this paper is to study the co-maximal hypergraph of dihedral groups, $Co_\mathcal{H}(D_n)$. We examine some of the structural properties viz., diameter, girth and chromatic number of $Co_\mathcal{H}(D_n)$. Also, we provide characterizations for  hypertrees, star structures and 3-uniform hypergraphs of  $Co_\mathcal{H}(D_n)$. Further,  we discuss the possibilities of $Co_\mathcal{H}(D_n)$ which can be embedded on the plane, torus and projective plane.
\end{abstract}
\section{Introduction}	
	Hypergraph is a generalization of graph, allowing the analysis of multiple relationships rather than just pair-wise relations. The notion of hypergraphs has been introduced by  C. Berge \cite{berge}. The study of hypergraphs on algebraic structures is an emerging area to extend some prominent results from graph theory. In \cite{cameron2021graphs}, P. J. Cameron introduced different types of graphs on groups whose edges reflect the group structures in some way. S. Akbari \textit{et.al.}  \cite{akbari} introduced the concept of co-maximal graph on subgroups of a group and they characterized all finite groups whose
	co-maximal graphs are connected. Later, in \cite{comaxzn}, A. Das \textit{et.al.} studied and characterized various properties like diameter, domination number, perfectness, hamiltonicity,
	etc. of the co-maximal graph on subgroups of cyclic groups. Recently, M. Saha and A. Das studied the co-maximal graph on subgroups of dihedral groups and  proved some of the isomorphism results related to it in \cite{codn}.
	
		A \textit{hypergraph} $\mathcal{H}$ is a pair $(V(\mathcal{H}),E(\mathcal{H}))$, where $V(\mathcal{H})$ is a set of  vertices and $E(\mathcal{H})$ is a set of hyperedges, where each hyperedge is a subset of $V(\mathcal{H})$. A hypergraph $\mathcal{H}'=(V'(\mathcal{H}'),E'(\mathcal{H}'))$   is called a \textit{subhypergraph} of  $\mathcal{H} = (V(\mathcal{H}),E(\mathcal{H}))$ if  $V'(\mathcal{H}')  \subseteq V(\mathcal{H})$ and $E'(\mathcal{H}') \subseteq E(\mathcal{H})$. A \textit{path in a hypergraph $\mathcal{H}$} is an alternating sequence of distinct
		vertices and edges of the form $v_1e_1v_2e_2...v_k$ such that $v_i,v_{i+1}$ is in $e_i$ for all $1 \leq i \leq k-1$. The \textit{cycle} 
	is a path whose first vertex is the same as the last vertex. The \textit{length} of a path is the
	number of hyperedges in the path. A hypergraph is said to be \textit{connected} if there exists a path between any two pair of vertices, otherwise 
	it is called a \textit{disconnected hypergraph}. The \textit{distance} between two vertices is the minimum length of the 
	path connecting these two vertices. The \textit{diameter} of a hypergraph is the maximum 
	distance among all pairs of vertices.  The \textit{girth} of a hypergraph is the length of a shortest cycle it contains. A hypergraph is called a \textit{star} if there is a vertex which belongs to all hyperedges.
	The \textit{incidence graph (or bipartite representation)} $\mathcal{I}(\mathcal{H})$ of $\mathcal{H}$ is a bipartite graph with vertex set $V(\mathcal{H}) \cup E(\mathcal{H})$ and a vertex $v \in V(\mathcal{H}) $ is adjacent to a vertex $e \in E(\mathcal{H})$ iff $v \in e$ in $\mathcal{H}$. A hypergraph $\mathcal{H}$ is called $r-$uniform, where $r$ is an integer, if for each edge
$e \in E(\mathcal{H}), |e| = r\, (r \geq 2)$. 	A \textit{proper vertex-coloring} (often simply called a proper coloring) of a hypergraph
$\mathcal{H}$ is an assignment of colors to the vertices of $\mathcal{H}$ such that no hyperedge contains all vertices of the
same color. The \textit{chromatic number} of $\mathcal{H}$, denoted by $\chi(\mathcal{H})$, is the minimum
number of colors needed for a proper vertex-coloring of $\mathcal{H}$. For more details on graphs and hypergraphs, one may refer \cite{voloshinhypergraphs,white1985graphs}, etc.

 In Section 2, we have introduced and studied the \textit{co-maximal hypergraph} $Co_\mathcal{H}(D_n)$  of dihedral groups and analyzed some of its structural properties viz., diameter, girth and chromatic number of $Co_\mathcal{H}(D_n)$. Also, we have characterized hypertrees, star hypergraphs and 3-uniform hypergraphs of  $Co_\mathcal{H}(D_n)$ in terms of $n$. We have obtained some results where there is a significant difference between some properties like girth, chromatic number, etc of co-maximal graph of $D_n$ and co-maximal hypergraph of $D_n$. In Section 3, we have discussed the possibilities of $Co_\mathcal{H}(D_n)$ which can be embedded on the plane, torus and projective plane.
 
 Note that dihedral group as a geometrical object is described by rotational and reflective symmetries and is an interesting class of non-abelian group. Hence, studying the co-maximal hypergraph of dihedral groups can add a different perspective.

\section{Co-maximal Hypergraph on $D_n$ and its structural properties}
	In this section, we introduce the concept of co-maximal hypergraphs of groups. Also, we analyze some of the structural properties, viz., diameter, girth and chromatic number of the co-maximal hypergraph  $Co_\mathcal{H}(D_n)$. Moreover, we characterize hypertrees, star hypergraphs and 3-uniform hypergraphs of  $Co_\mathcal{H}(D_n)$. In \cite{codn}, A. Das and M. Saha introduced the co-maximal subgroup graph $\Gamma(G)$ of a group $G$ as follows:
	\begin{defn}\cite{codn}
		Let $G$ be a group and $S$ be the collection of all non-trivial proper subgroups of $G$. The \textit{co-maximal subgroup graph} $ \Gamma(G)$ of a group $G$ is defined to be a graph with $S$ as the set of vertices and two distinct vertices $H$ and $K$ are adjacent if and only if $HK=G$. The \textit{deleted co-maximal subgroup graph of $G$}, denoted by $\Gamma^\star(G)$, is defined as the graph obtained by removing the isolated vertices from $\Gamma(G)$.
	\end{defn}
	\noindent Motivated from \cite{cameron2023hypergraphs}, we have defined a hypergraph as follows:
	\begin{defn}
		Let $G$ be a group and $S$ be the set of all non-trivial proper subgroups of $G$. The \textit{co-maximal hypergraph of $G$}, denoted by $Co_{\mathcal{H}}(G)$, is a hypergraph whose vertex set, $V=\{H \in S \, | \, HK=G \,\, \text{for some} \, K \in S\}$ and $E \subseteq V$ is a hyperedge if and only if \begin{enumerate}
			\item for distinct $H,K \, \in E$, $HK = G.$
			\item there does not exist  $E' \supset E$ which  satisfies (1).
		\end{enumerate}
	\end{defn}

	\begin{exmp}

			Consider the Klein-4 group, $$V_4=\{e,a,b,c \, | a^2=b^2=c^2=e, \,  ab=c=ba, \, ac=b=ca, \, bc=a=cb\}.$$ Then, the vertex set of $Co_\mathcal{H}(V_4)$ is $V= \{ \{e,a\}, \{e,b\}, \{e,c\}\}$ and the hyperedge set is $\{ \{\{e,a\}, \{e,b\}, \{e,c\}\}\}$. 
		
	\begin{figure}[h]
		\centering\subfloat[$Co_\mathcal{H}(V_4)$]{	\begin{tikzpicture}[scale=.3]
				\draw[thick] (0,0) ellipse (7cm and 2cm);
				
				\filldraw (-1cm, 0.5cm) circle (2pt) node[right] {$\{e,b\}$};    
				\filldraw (-2.0cm, -0.3cm) circle (2pt) node[ left] {$\{e,a\}$}; 
				\filldraw (2.5cm, -0.3cm) circle (2pt) node[right] {$\{e,c\}$};  
		\end{tikzpicture}}
		\label{G}
		\hspace{.50 cm}
		\subfloat[$\Gamma^\star(V_4)$]{
			\begin{tikzpicture}[scale=1.5, every node/.style={circle, fill=black, minimum size=6pt, inner sep=0pt}]
				
				\node[label=left:{$\{e,a\}$}] (u1) at (-1,1) {};
				\node[label=below:{$\{e,b\}$}] (u2) at (0,0) {};
				\node[label=right:{$\{e,c\}$}] (u3) at (1,1) {};

				\draw (u1) -- (u2);
				\draw (u1) -- (u3);
				\draw (u2) -- (u3);

		\end{tikzpicture}}
		\caption{}
		\label{klein4group}

	\end{figure}
		
		\end{exmp}
\begin{note}
$V_4 \cong \mathbb{Z}_2 \times \mathbb{Z}_2 \cong D_2$.
\end{note}
\begin{exmp}
	Consider the dihedral group $D_4$ of order 8, $D_4= <a,b\,|a^4=e=b^2, bab^{-1}=a^{-1}>$.

		\noindent The vertex set of $Co_\mathcal{H}(D_4)$ is $ V=\{H_3,H_4,H_5, H_6,H_7,H_8,H_9\}$, where 
	 $ H_3=<b>, H_4=<ab>,$ $H_5=<a^2b>, H_6=<a^3b>, H_7=<a>, H_8=<a^2,ab>,$  $ H_9=<a^2,b> $. The hyperedge set of 	$Co_\mathcal{H}(D_4)$	is $\{e_1,e_2,e_3,e_4,e_5\}$, where  
	$e_1 = \{H_3,H_7,H_8\}$,
	$e_2=\{H_4,H_7,H_9\}$,
	$e_3=\{H_5,H_7,H_8\}$, $e_4=\{H_6,H_7,H_9\}$ and $e_5=\{H_7,H_8,H_9\}$.

	\begin{figure}[h]
	\centering\subfloat[$Co_\mathcal{H}(D_4)$]{	\begin{tikzpicture}[scale=1]
			
			\draw[rotate=30,thick] (.8,-0.5) ellipse (.7 and 2);
			\draw[thick] (1,0) ellipse (3 and .5);
			\draw[thick] (1,-.8) circle (1.5cm);
			\draw[thick] (2,.2) circle (1.2cm);
			\draw[thick] (2,-.5) circle (1.2cm);

			\node at (-2.2, 0.2) {$e_1$};
			\node at (-.5, 1.5) {$e_2$};
			\node at (2.5, 1.5) {$e_3$};
			\node at (-.6, -1.5) {$e_4$};
			\node at (3, -1.5) {$e_5$};


			\filldraw (-1,-0.2cm) circle (2pt) node[right] {${H_3}$};
			\filldraw (2,.9cm) circle (2pt) node[above] {${H_5}$};
			\filldraw (1.1, -0.2cm) circle (2pt) node[right] {${H_7}$};
			\filldraw (0.2, 1cm) circle (2pt) node[above] {${H_4}$};
			\filldraw (0, -1.2cm) circle (2pt) node[below] {${H_6}$};
			\filldraw (2, -1.2cm) circle (2pt) node[left] {${H_9}$};
			\filldraw (3, -0.2cm) circle (2pt) node[left] {${H_8}$};

	\end{tikzpicture}	}
	\label{G}
	\hspace{.50 cm}
	\subfloat[$\Gamma^\star(D_4)$]{
		\begin{tikzpicture}[scale=1.5, every node/.style={circle, fill=black, minimum size=6pt, inner sep=0pt}]
			
			\node[label=above:$H_4$] (u1) at (-.5,2) {};
			\node[label=above:$H_6$] (u2) at (.5,2) {};
			\node[label=left:$H_9$] (u3) at (0,0) {};
			
			\node[label=above:$H_3$] (v1) at (2.5,2) {};
			\node[label=above:$H_5$] (v2) at (1.5,2) {};
			\node[label=right:$H_8$] (v3) at (2,0) {};
			\node[label=below:$H_7$] (v4) at (1,1) {};
			
			\draw (u3) -- (v3);
			\draw (v4) -- (v3);
			\draw (v4) -- (u3);
			\draw (u3) -- (u1);
			\draw (v4) -- (u1);
			\draw (u3) -- (u2);
			\draw (v4) -- (u2);
			\draw (v3) -- (v2);
			\draw (v4) -- (v2);
			\draw (v3) -- (v1);
			\draw (v4) -- (v1);
			
	\end{tikzpicture}}
	\caption{}
	\label{d41}

\end{figure}

\end{exmp}
\begin{flushright}
	$\square$
\end{flushright}
	\begin{remark}\label{remark1}
	$Co_{\mathcal{H}}(G)$ is the clique hypergraph of $\Gamma^\star(G)$, i.e., the hyperedges of $Co_{\mathcal{H}}(G)$  are the maximal cliques of $\Gamma^\star(G)$. 
\end{remark}
 For a positive integer $n \geq 1$, the \textit{dihedral group of order $2n$} is denoted by $D_n$ and  is defined as
$$D_n=<a,b\, | a^n=e, \, b^2=e,\,  bab^{-1}=a^{-1}>.$$ 

\noindent The following theorem explicitly describes the subgroups of dihedral groups. 
 \begin{theorem}\cite{conrad2009dihedral}\label{Kconrad1}
	Every subgroup of $D_n$ is cyclic or dihedral. A complete listing of the subgroups is as follows:
	\begin{enumerate}
		\item $< a^r > $ with index 2$ r$, where $r \mid n$.
		\item $ <a^r, a^i b> $ with index $r$, where $r \mid n $ and $0 \leq i \leq r - 1$.
	\end{enumerate}
	Every subgroup of $D_n$ occurs exactly once in this listing.
\end{theorem}
\begin{remark}
	\begin{enumerate}
		\item A subgroup of $D_n$ is said to be of \textbf{Type (1)} if it is cyclic as stated in 1 of Theorem \ref{Kconrad1}.
		\item A vertex of $Co_\mathcal{H}(D_n)$ is said to be of \textbf{Type (1)} if it is a subgroup of $D_n$ of  \textbf{Type (1)}.
		\item A subgroup of $D_n$ is said to be of \textbf{Type (2)} if it is dihedral subgroup as stated in 2 of Theorem \ref{Kconrad1}.
		\item A vertex of $Co_\mathcal{H}(D_n)$ is said to be of \textbf{Type (2)} if it is a subgroup of $D_n$ of  \textbf{Type (2)}.
	\end{enumerate}
\end{remark}

\begin{remark}\label{remarkcard}
	The following observations from \cite{mitkari} are useful for the subsequent results.  Here, for subgroups $H,K$ of $D_n$, $H \vee K=$ $<H \cup K>$ and $H \wedge K = H \cap K$.
		\begin{enumerate}
		\item Let $H=<a^{n_1}>$ and $K=<a^{n_2}>$, where $\mid H\mid =m_1=\frac{n}{n_1}, \, \mid K \mid =m_2=\frac{n}{n_2}$, then $H \vee K=$ $<a^{(n_1,n_2)}>$ and $H\wedge K=<a^{[n_1,n_2]}>$, where $\mid H \vee K \mid =[m_1,m_2]=\frac{n}{(n_1,n_2)},$ $ \mid H \wedge K \mid =(m_1,m_2) = \frac{n}{[n_1,n_2]}$.
		\item  Let $H=<a^{n_1}>$ and $K=<a^{n_2},a^ib>$, where $\mid H \mid = m_1=\frac{n}{n_1}, \mid K \mid =m_2=\frac{2n}{n_2}$, then $H \vee K = <a^{(n_1,n_2)},a^ib>$ and $H \wedge K=<a^{[n_1,n_2]}>$, where $\mid H \vee K \mid = [m_1,m_2]=\frac{2n}{(n_1,n_2)}$, $\mid H \wedge K \mid = (m_1,m_2)=\frac{n}{[n_1,n_2]}$.
		\item Let $H=<a^{n_1},a^ib>$ and $K=<a^{n_2},a^jb>$, where $\mid H \mid=m_1=\frac{2n}{n_1}, \mid K \mid=m_2=\frac{2n}{n_2}$, then $H \vee K =<a^{(n_1,n_2)},a^ib>$, where $\mid H \vee K \mid=[m_1,m_2]=\frac{2n}{(n_1,n_2)}$ and, \begin{enumerate}
			\item If $n_1x +n_2y=i-j$ has no integer solution, then $H\wedge K = <a^{[n_1,n_2]}>$, where $\mid H \wedge K \mid =  (m_1,m_2)=\frac{n}{[n_1,n_2]}$.
			\item If $n_1x +n_2y=i-j$ has an integer solution, then $H\wedge K = <a^{[n_1,n_2]},a^{i-n_1x_0}b>$, where $\mid H \wedge K \mid =  (m_1,m_2)=\frac{2n}{[n_1,n_2]}$ and $ (x_0,y_0)$ is an integer solution of the equation $n_1x +n_2y=i-j$.
		\end{enumerate}
		
	\end{enumerate}
\end{remark}
	\begin{note}
	Since $Co_\mathcal{H}(D_n)$ is empty for $n=1$, we exclude this case and consider $n \geq 2$ throughout the article.
\end{note}
\begin{prop}\cite{codn} \label{stargraph}
	$\Gamma^{\star}(D_n)$ is a star if and only if $n$ is an odd prime power.
\end{prop}
\begin{remark} \label{oddprime}
	From the above proposition, we have $\Gamma^{\star}(D_n) \cong Co_\mathcal{H}(D_n)$ for an odd prime power $n$.
\end{remark}

The following theorem follows from Theorem 4 of \cite{codn}. But for the sake of completeness we are giving the proof which we use in the subsequent results.
\begin{theorem}\label{vertexset}
	 $Co_\mathcal{H}(D_n)$ is non-empty. Moreover, 	all the non-trivial proper subgroups of $D_n$ constitutes the vertex set of $Co_\mathcal{H}(D_n)$ if and only if $n$ is a square-free.
	
\end{theorem}

\begin{proof}
	Consider the subgroup $<a>$ of $D_n$. Observe that $<a>\cdotp<b>=D_n$. Hence, $<a>,\\<b> \in V(Co_\mathcal{H}(D_n))$ and therefore, $Co_\mathcal{H}(D_n)$ is non-empty.
	
	 All the  \textbf{Type (2)} subgroups of $D_n$ are in the vertex set of $Co_\mathcal{H}(D_n)$ as their product with $<a>$ is equal to $D_n$. 
	 Also, as \textbf{Type (1)} subgroups of $D_n$ are normal, by Remark \ref{remarkcard}.1, $\mid <a^{r_1}> \cdotp <a^{r_2}> \mid =$ $ \mid <a^{r_1}> \vee <a^{r_2}> \mid = \frac{n}{(r_1,r_2)}$. Thus, the product of any two \textbf{Type (1)} subgroups of $D_n$ is not equal $D_n$. Now, for the subgroup $<a^{r}>$ of $D_n$, consider the following cases: \\
	  \textbf{Case 1.} Suppose all prime divisors of $n$ are divisors of $r$. Then,  by Remark \ref{remarkcard}.2,\\ $\mid <a^{r}> \cdotp <a^{r_2}, b> \mid = \mid <a^{r}> \vee <a^{r_2}, b> \mid =  \mid <a^{(r,r_2)},b> \mid $ = $\frac{2n}{(r,r_2)} \neq 2n$ and therefore, $<a^r>$ does not belong to $ V(Co_\mathcal{H}(D_n))$. \\
	  \textbf{Case 2}. Suppose $p_1$ is a prime divisor of $n$ which is not a divisor of $r$. By Remark \ref{remarkcard}.2, \\$\mid <a^{r}> \cdotp <a^{p_1}, b> \mid = \mid <a^{r}> \vee <a^{p_1}, b> \mid =  \mid <a^{(r,p_1)},b> \mid $ = $\frac{2n}{(r,p_1)} = 2n$ and therefore, $<a^r>$ belongs to $ V(Co_\mathcal{H}(D_n))$. \\
	 Thus,  $<a^r>$ is not in the $V(Co_\mathcal{H}(D_n))$ if and only if all prime divisors of $n$ are divisors of $r$. Therefore, from the above cases we can conclude that all the non-trivial proper subgroups of $D_n$ belongs to  $V(Co_\mathcal{H}(D_n))$ iff $n$ is a square-free.
\end{proof}
\begin{remark}\label{powerof2} 
	 For $n=2^\alpha,$ where $\alpha$ is a positive integer greater than or equal to 1, the only \textbf{Type (1)} vertex of $Co_{\mathcal{H}}(D_n)$ is $<a>$ and $<a>$ belongs to all the hyperedges of $Co_{\mathcal{H}}(D_n)$. \\
	 For $n=2$, $Co_\mathcal{H}(D_n)$ consists of a single hyperedge containing the vertices $<a>,<b>$ and $<ab>$. \\
	 For $n=2^\alpha,$ where $\alpha$ is a positive integer greater than 1, if  $r_1,r_2$ are integers greater than $1$ such that $r_1 \neq n, r_2 \neq n$,  $r_1 \mid n$, $r_2 \mid n$ and $r_1 \mid r_2$, then the following results hold:
	\begin{enumerate}
		\item  For the vertices $<a^{r_1},a^ib>$ and $<a^{r_2},a^jb>$ where  $0\leq i \leq r_1 -1$ and $0\leq j \leq r_2 -1$,   $<a^{r_1},a^ib> \cap <a^{r_2},a^jb>$ is either $<a^{r_2}>$ or $<a^{r_2},a^jb>$.
		\item $ <a^2,b> \cdotp <a^{r_1}, a^ib> = D_n$ and $ <a^2,ab> \cdotp <a^{r_1}, a^ib> \neq  D_n$ for $i \equiv 1 (mod \, 2)$.
		\item  $ <a^2,ab> \cdotp <a^{r_1}, a^jb> = D_n$ and $ <a^2,b> \cdotp <a^{r_1}, a^jb> \neq D_n$ for $j \equiv 0 (mod \, 2)$.
		\item If  $r_1\neq 2,  r_2 \neq 2$ and  $r_1 \mid r_2$, then  $<a^{r_1},a^ib>\cdotp <a^{r_2},a^jb> \neq D_n$, where $0 \leq i \leq r_1 -1, \, \, 0 \leq j \leq r_2 -1$. 
	\end{enumerate}

 Thus,  for $n=2^\alpha,$ where $\alpha$ is a positive integer greater than 1, we get the structure of hyperedges of $Co_\mathcal{H}(D_n)$  as follows: \begin{enumerate}
\item $e=\{<a>,<a^2,b>,<a^2,ab>\}$,
\item $e^{'}=\{<a>,<a^2,b>,<a^r,a^jb>\}$, where  $j \equiv 1 (mod \, 2)$,
\item $e^{''}=\{<a>,<a^2,ab>,<a^r,a^jb>\}$, where $j \equiv 0 (mod \, 2)$.
\end{enumerate}
where $r$ is a non-negative integer dividing $n$.
\end{remark}
\begin{prop}\label{diam}
	If $\mathcal{G}$ is connected, then $Cl_\mathcal{H}(\mathcal{G})$ is also connected. Moreover, $diam(\mathcal{G})= diam(Cl_\mathcal{H}(\mathcal{G}))$.
\end{prop}
\begin{proof}
	For a graph $\mathcal{G}$,  observe that $ \mid V(Cl_\mathcal{H}(\mathcal{G}))\mid = \mid V(\mathcal{G}) \mid -$ number of isolated vertices in $\mathcal{G}$.
	Let $x,y$ be any two vertices in $\mathcal{G}$ with $dist(x,y)=k$ in $\mathcal{G}$. Then there exists a shortest path $P$ from $x$ to $y$ of length $k$ in $\mathcal{G}$. Let $P = x \, e_1 \, a_1 \,  e_2 \, a_2 \, \cdots \, a_{k-1}\, e_k \,y$. Observe that no non-consecutive vertices in the path belongs to a clique. Hence, $dist(x,y)=k$ in $Cl_\mathcal{H}(\mathcal{G})$.
\end{proof}
The following theorem follows directly from the above proposition and from Theorem 6 of \cite{codn}.
\begin{theorem}
	$Co_{\mathcal{H}}(D_n)$ is connected with $diam(Co_{\mathcal{H}}(D_n)) \leq 3$. In particular, 	$$
	diam(Co_\mathcal{H}(D_n)) = \begin{cases*} 1 \,\, \text{if $n =2$,}\\
		2 \,\,  \, \text{if} \, n= p^{\alpha}, \, \text{where $p$ is a prime,  $\alpha \geq 1$ and $n \neq 2$, }\\
		3 \,\, \text{otherwise}.
	\end{cases*}$$
	
\end{theorem}

 In \cite{codn}, A. Das and M. Saha established that  $\Gamma^{\star}(D_n)$ is a star if and only if $n$ is an odd prime power if and only if $\Gamma^{\star}(D_n)$ is a tree. But, in the case of co-maximal hypergraph on $D_n$, we have proved that $Co_\mathcal{H}(D_n)$ is a star hypergraph if and only if $n$ is a power of a prime if and only if  $Co_\mathcal{H}(D_n)$ is a hypertree.\\
 Note that every star hypergraph is a hypertree. For characterizing  $Co_{\mathcal{H}}(D_n)$ that are hypertrees and star hypergraphs, we need the following definitions and results.
	\begin{defn}\cite{voloshinhypergraphs}
		A \textit{host graph} for a hypergraph is a connected graph $G$ on the same vertex set such that
		every hyperedge induces a connected subgraph of $G$. A hypergraph $\mathcal{H} = (X,\mathcal{D})$ is called a \textit{hypertree} if there exists a host tree
		$T = (X,E)$ such that each edge $D \in \mathcal{D}$  induces a subtree in $T$. 
	\end{defn}
	\begin{defn}\cite{voloshinhypergraphs}
		A hypergraph $\mathcal{H}$ has the \textit{Helly property (is Helly, for short)} if for every subfamily of its edges the following implication holds:\\ If every two edges of the subfamily have a  non-empty intersection, then the whole subfamily has a non-empty intersection.
	\end{defn}

	\begin{lemma}\cite{voloshinhypergraphs}\label{hellyhypertree}
		Every hypertree is a Helly hypergraph.
	\end{lemma}

\begin{theorem}\label{starhypergraph}
	For $Co_\mathcal{H}(D_n)$, the following statements are equivalent:
	\begin{enumerate}
		\item $Co_{\mathcal{H}}(D_n)$ is a hypertree.
		\item $n=p^\alpha$, where $p$ is a prime and $\alpha \geq 1$.
		\item $Co_{\mathcal{H}}(D_n)$ is a star hypergraph.
	\end{enumerate}
\end{theorem}
\begin{proof}
	$1 \Rightarrow 2$. Assume that   $Co_{\mathcal{H}}(D_n)$ is a hypertree and  $n$ is not a power of a prime. Consider the following cases:\\
	\textbf{Case 1.} Suppose that $n$ is even. Consider set  $S =\{e_1.e_2,e_3,e_4\}$ of hyperedges of   $Co_{\mathcal{H}}(D_n)$, where $e_1=\{<a>,<a^2,b>,<ab>\}$, $e_2=\{<a>,<a^2,ab><b>\}$, \\$e_3 \supseteq \{<a>,<a^{p_1},b>,<a^2,b>,<a^2,ab>\} =\alpha$, $e_4 \supseteq \{<a^{p_1}>,<a^2,b>,<a^2,ab>\}=\beta$ and $p_1$ is an odd prime divisor of $n$. Observe that any two hyperedges of $S$ have a non-empty intersection but as $e_1 \cap e_2 \cap e_3 \cap e_4 = e_1 \cap e_2 \cap \alpha \cap \beta$ = $\emptyset$,  $S$ does not satisfy the Helly property. Thus, by Lemma \ref{hellyhypertree},  $Co_{\mathcal{H}}(D_n)$ is not a hypertree, which is a contradiction.
	
	\noindent \textbf{Case 2.} Suppose that $n$ is odd and $n=p_1^{\alpha_1}p_2^{\alpha_2}\ldots p_k^{\alpha_k}$ where $p_i$'s are prime divisors of $n$ and $\alpha_i$'s are non-negative integers. Consider the set $S=\{e_1,e_2,e_3\}$ of hyperedges of  $Co_{\mathcal{H}}(D_n)$, where $e_1 = \{<a>,<a^{p_1^{\alpha_1}p_2^{\alpha_2}},b>, <a^{p_3^{\alpha_3}\ldots p_{k-1}^{\alpha_{k-1}}},b>, <a^{p_k},b>\}$, \\
	$e_2 = \{<a>, <a^{p_1^{\alpha_1}},b>, <a^{p_2^{\alpha_2}},b> <a^{p_3^{\alpha_3} \ldots p_{k-1}^{\alpha_{k-1}} p_k^{\alpha_k}},b>\}$, and\\$e_3= \{<a^{p_1}>,<a^{p_2^{\alpha_2}},b>,<a^{p_3^{\alpha_3}\ldots p_{k-1}^{\alpha_{k-1}}},b><a^{p_k},b>\}$. Observe that any two hyperedges of $S$ have a non-empty intersection but as $e_1 \cap e_2 \cap e_3 =\emptyset$, $S$ does not satisfy Helly property. Thus, by Lemma \ref{hellyhypertree}, $Co_{\mathcal{H}}(D_n)$ is not a hypertree, which is a contradiction.\\
		$2 \Rightarrow 3$.
			Assume that $n=p^\alpha$ where $p$ is a prime and $\alpha \geq 1$.\\
		If $n=2$, then $Co_{\mathcal{H}}(D_n)$ is a hypergraph consisting of a single hyperedge and hence a star hypergraph. \\
		If $n=p^{\alpha}$, where $p$ is a prime, $ \alpha \geq 2$ and $n\neq 2$, then by Remark \ref{oddprime} and Remark \ref{powerof2} $Co_{\mathcal{H}}(D_n)$ is a star hypergraph.\\
			$3 \Rightarrow 1$. This implication is straightforward.
\end{proof}
In \cite{codn}, A. Das and M. Saha proved that the girth of $\Gamma(D_n)$ is 3 for all $n \geq 3$ except odd prime powers. Note that $\mathcal{G}$ has a cycle does not imply that $Cl_\mathcal{H}(\mathcal{G})$ has a cycle. In the following result, we have proved that the girth of co-maximal hypergraph of $D_n$ is either 2 or $\infty$.
\begin{theorem}
	The girth $gr(Co_{\mathcal{H}}(D_n))$ of $Co_{\mathcal{H}}(D_n)$ is either 2 or  $\infty$. In particular, $$
	gr(Co_{\mathcal{H}}(D_n))  = \begin{cases*} \infty,  \,\, \text{if $n=2$ or $n =p^\alpha$, where $p$ is an odd prime and $\alpha$ is a positive integer,}\\
		2,  \,\, \text{otherwise}.
	\end{cases*}$$
\end{theorem}
\begin{proof}
	We consider the following cases:\\
	\noindent $\underline{\textbf{Case 1.}}$ If $n=2$, then  by Remark \ref{powerof2}, $gr(Co_{\mathcal{H}}(D_n)) =\infty$.

		\noindent $\underline{\textbf{Case 2.}}$ If $n=p^\alpha$, where $p$ is an odd prime and $\alpha$ is a positive integer, then by Remark \ref{oddprime},  $gr(Co_{\mathcal{H}}(D_n))=\infty$.
	
	\noindent $\underline{\textbf{Case 3.}}$
If $n=2^\alpha$, where $\alpha \geq 2$, then, by Remark \ref{powerof2}, there exist two distinct hyperedges $e_1$ and $e_2$ such that $e_1$ contains $<a>,$ $<a^2,b>$ and  $<a^2,ab>$, and  $e_2$ contains $<a>,<a^2,b>$ and  $<ab>$. Therefore, $<a> e_1 <a^2,b> e_2 <a>$ is a shortest cycle of length 2 and consequently, $gr(Co_{\mathcal{H}}(D_n))=2$.
		\begin{figure}[h]
	\centering\subfloat[]{	\begin{tikzpicture}[scale=.7]
			
			\draw[thick] (-1.7,.8) ellipse (3 and 1);
			
			\draw[thick] (-1,0) ellipse (1.5 and 2.5);
			
			\node at (-3.2, 2) {$e_1$};
			\node at (-1, 2.7) {$e_2$};


			\filldraw (-1.6, 1.2cm) circle (2pt) node[right] {${<a>}$};
			\filldraw (-1.6, .5cm) circle (2pt) node[right] {${<a^2,b>}$};
			\filldraw (-1, -1cm) circle (2pt) node[below] {${<ab>}$};
			\filldraw (-3.4, 1.2cm) circle (2pt) node[below] {${<a^2,ab>}$};
	\end{tikzpicture}	}
	\label{G}
	\hspace{.50 cm}
	\subfloat[]{
		\begin{tikzpicture}[scale=.7]
			
			\draw[thick] (-1.7,.8) ellipse (3 and 1.5);
			
			\draw[thick] (-1,0) ellipse (1.5 and 3);
			
			\node at (-3.2, 2.3) {$e_1$};
			\node at (-1, 3.2) {$e_2$};


			\filldraw (-1.6, 1.2cm) circle (2pt) node[right] {${<a>}$};
			\filldraw (-1.6, .5cm) circle (2pt) node[right] {${<a^{p_1},b>}$};
			\filldraw (-1, -1.3cm) circle (2pt) node[below] {${<a^{p_2},ab>}$};
			\filldraw (-3.4, 1.2cm) circle (2pt) node[below] {${<a^{p_2},b>}$};
	\end{tikzpicture}}
	\caption{}
	\label{girth2}

\end{figure}

	\noindent $\underline{\textbf{Case 4.}}$ 	
		If $n$ is not a power of a prime, then there exist atleast two distinct prime divisors, say $p_1,p_2$ of $n$. Without loss of generality, assume $p_2 \neq 2$. Then, for the vertices $<a>$ and $<a^{p_1},b>$, there exist two distinct hyperedges $e_1$ and $e_2$ such that $e_1$ contains $<a>,<a^{p_1},b>$ and  $<a^{p_2},b>$, and  $e_2$ contains $<a>,<a^{p_1},b>$ and $<a^{p_2},ab>$. Therefore, $<a> e_1 <a^{p_1},b> e_2 <a>$ is a shortest cycle of length 2 and consequently, $gr(Co_{\mathcal{H}}(D_n))=2$. 
\end{proof}
$\chi(G)$ and $\chi(Cl_\mathcal{H}(G))$ need not be equal in general. But for any graph $\mathcal{G}$,  $\chi(Cl_\mathcal{H}(\mathcal{G}))$ and $\chi_c(\mathcal{G})$ coincides, where $\chi_c(\mathcal{G})$ denotes the clique chromatic number of a graph $\mathcal{G}$, i.e., the smallest number of colors in a vertex coloring so that no maximal clique is monochromatic.
\noindent In \cite{codn}, it is proved that the chromatic number of $\Gamma(D_n)$ is as follows:
$$
\chi(\Gamma(D_n))  = \begin{cases*} \pi(n) + 1,  \,\, \text{if $n$ is odd,}\\
	\pi(n) +2,  \,\, \text{if $n$ is even}.
\end{cases*}$$
 In the following result, we have established  that the chromatic number of co-maximal hypergraph  $Co_{\mathcal{H}}(D_n)$ of $D_n$, for all $n\geq 2$, is 2. Thus, $\chi(Co_{\mathcal{H}}(D_n))$ and $\chi(\Gamma(D_n))$ coincides only when $n$ is a power of an odd prime as in this case $Co_\mathcal{H}(D_n) \cong \Gamma(D_n)$ and in all other cases the chromatic number differs.
\begin{theorem}
		The chromatic number,
	$\chi(Co_{\mathcal{H}}(D_n)) = 2.$
\end{theorem}

\begin{proof}
		Divide the vertex set of $Co_{\mathcal{H}}(D_n)$ into two sets $A \, \text{and} \, B$, where $A$ is the set of all \textbf{Type (1)} vertices and $B$ is the set of all \textbf{Type (2)} vertices. 
	Let $e_1$ be a hyperedge of 	$Co_{\mathcal{H}}(D_n)$.  Since no two product of \textbf{Type (1)} vertices is equal to $D_n$, $e_1$ cannot contain all vertices from $A$ only. Also note that the product of the vertex $<a>$ with any \textbf{Type (2)} vertex is equal to $D_n$. So, $e_1$  cannot contain all vertices from $B$ only because of  the maximality of $e_1$. 
		 Hence, each hyperedge has at least one vertex from both A and B. Assign the color $c_1$ to the vertices in $A$ and the color $c_2$ to the vertices in $B$. This is a proper coloring of $Co_{\mathcal{H}}(D_n)$. Consequently,  
		$\chi(Co_{\mathcal{H}}(D_n)) = 2.$
\end{proof}

As $\Gamma^{\star}(D_n)$ is a simple graph, i.e., a 2-uniform hypergraph,  it is interesting to know when  $Co_{\mathcal{H}}(G)$ on a group $G$ is a $k$-uniform hypergraph. In the following result, we have settled this question for $k=3$ and $G=D_n$.

\begin{theorem}\label{uniformhyp}
	$Co_{\mathcal{H}}(D_n)$ is a 3-uniform hypergraph if and only if $n=2^\alpha$, where $\alpha$ is a positive integer.
\end{theorem}
\begin{proof}
	If $n=2^\alpha$, where $\alpha$ is a positive integer, then by Remark \ref{powerof2}, $Co_{\mathcal{H}}(D_n)$ is a 3-uniform hypergraph. 
	
Conversely, suppose that 	$Co_{\mathcal{H}}(D_n)$ is a 3-uniform hypergraph. If $n$ is an odd prime power, then	$Co_{\mathcal{H}}(D_n)$ is a 2-uniform hypergraph, but not a 3-uniform.  Now, assume that $n$ is not a power of two. Consider the following cases:\\
\noindent\textbf{Case 1.} If $n$ is even, then there exists a hyperedge $e_1 = \{<a>,<a^2,b>,<ab>\}$ and another hyperedge $e_2 \supseteq \{<a>,<a^2,b>,<a^2,ab>,<a^p,b>\}$, where $p$ is an odd prime divisor of $n$. Thus, $\mid e_1 \mid = 3$ and $\mid e_2 \mid \geq 4$ and consequently,	$Co_{\mathcal{H}}(D_n)$ is not a 3-uniform hypergraph, which is a contradiction.\\
\noindent\textbf{Case 2.} If $n$ is odd, then $e_1 =\{<a>,<b>\}$ is a hyperedge of $Co_{\mathcal{H}}(D_n)$ and there exists a hyperedge $e_2 \supseteq \{<a>,<a^p,b>,<a^q,b>\}$ of	$Co_{\mathcal{H}}(D_n)$, where $p$ and $q$ are distinct  prime divisors of $n$. Thus, $\mid e_1 \mid =2$ and $\mid e_2 \mid \geq 3$ and consequently, $Co_{\mathcal{H}}(D_n)$ is not a 3-uniform hypergraph, which is a contradiction.

\end{proof}

\section{Embedding of $Co_{\mathcal{H}}(D_n)$}

 An \textit{embedding of a graph} on a surface is a continuous and  one to one
function from a topological representation of the graph into the surface. We denote by $S_n$ the surface obtained from the sphere $S_0$  by adding $n$ handles. The number $n$ is called the \textit{genus of the surface} $S_n, n\geq 0$.  The \textit{orientable genus} of a graph $\mathcal{G}$, denoted by $g(\mathcal{G})$, is the minimum genus of
a surface in which $\mathcal{G}$ can be embedded.  A \textit{cross-cap} is a topological object formed by identifying opposite points on the boundary of a circle (or a disk) and is equivalent to gluing a Möbius strip into a hole in a surface. A surface obtained by adding $k$ crosscaps to $S_0$ 
is known as the non-orientable surface and we denote it by $N_k$. The number $k$ is called the crosscap of $N_k$. The non-orientable genus of
a graph $\mathcal{G}$, denoted by $\tilde{g}(\mathcal{G})$, is the smallest integer $k$ such that $\mathcal{G}$ can be embedded on $N_k$.  A graph is said to be \textit{planar} if it can be drawn on the plane in such a way that no edges intersect, except at a common end vertex.  A graph is said to be \textit{toroidal} if it  can be embedded on a torus and is called \textit{projective} if it can be embedded on a projective plane.   Further, note that if $H$ is a subgraph of a graph $\mathcal{G}$, then $g(H) \leq g(\mathcal{G})$ and $\tilde{g}(H) \leq \tilde{g}(\mathcal{G})$. A hypergraph is \textit{toroidal} if its incidence graph is toroidal and is  \textit{projective} if its incidence graph is projective. \\
  We are interested to study  $Co_{\mathcal{H}}(D_n)$ which can be embedded on plane, torus, projective plane, etc.  First, we discuss the planarity of $Co_{\mathcal{H}}(D_n)$. To analyze the planarity of $Co_{\mathcal{H}}(D_n)$, we need the following results.

\begin{theorem}\cite{walsh1975hypermaps}
	A graph $G$ is planar iff it contains no subdivision of $K_5$ or $K_{3,3}$.
\end{theorem}
\begin{theorem}\cite{walsh1975hypermaps} 
	A hypergraph is planar iff its incidence graph is planar.
\end{theorem}

\begin{theorem}\label{planarprimepower}
	$Co_{\mathcal{H}}(D_n)$  is planar if and only if $n=p^\alpha$, where $p$ is a prime and $\alpha$ is a positive integer.
\end{theorem}
\begin{proof}

	Suppose that $n$ is not a power of a prime. Then we will prove that $Co_{\mathcal{H}}(D_n)$  contains either $K_{3,3}$ or a subdivision of $K_{3,3}$. For proving this, consider the following cases: \\
	\noindent \textbf{Case 1.} Suppose that $n$ is even. For the vertices,  $H_1=<a>,H_2=<a^2,b>,H_3=<a^2,ab>,$ $K_1=<a^{p_1},b>, K_2=<a^{p_1},ab>$ and $K_3=<a^{p_1},a^2b>$ where $p_1$ is an odd prime divisor of $n$,   there exist three distinct hyperedges $e_1,e_2,e_3$ of $Co_{\mathcal{H}}(D_n)$ such that $e_1$ contains $H_1,H_2,H_3$ and $K_1$, $e_2$ contains $H_1,H_2,H_3$ and $K_2$, and $e_3$ contains $H_1,H_2,H_3$ and $K_3$. Therefore,  $\mathcal{I}(Co_{\mathcal{H}}(D_n))$ contains $K_{3,3}$. Hence, $Co_{\mathcal{H}}(D_n)$ is non-planar.
		
	\noindent\textbf{Case 2.} Suppose that $n$ is odd.\\
	\noindent\textbf{Subcase 2.1.}  Suppose that $\omega(n)=2$, i.e., $n=p_1^{\alpha_1}p_2^{\alpha_2}$, where $p_1,p_2$ are odd primes and $\alpha_1, \alpha_2$ are positive integers. Consider the following hyperedges of $Co_{\mathcal{H}}(D_n)$:  \\ $e_1=\{<a>,<a^{p_1^{\alpha_1}},a^{i_1}b>,<a^{p_2^{\alpha_2}},a^{j_1}b>\},  e_2=\{<a>,<a^{p_1^{\alpha_1}},a^{i_1}b>,<a^{p_2^{\alpha_2}},a^{j_2}b>\},$ \\ $ e_3=\{<a>,<a^{p_1^{\alpha_1}},a^{i_1}b>,<a^{p_2^{\alpha_2}},a^{j_3}b>\}, e_4=\{<a>,<a^{p_1^{\alpha_1}},a^{i_2}b>,<a^{p_2^{\alpha_2}},a^{j_2}b>\},$ \\  $ e_5=\{<a>,<a^{p_1^{\alpha_1}},a^{i_2}b>,<a^{p_2^{\alpha_2}},a^{j_1}b>\},$ $e_6=\{<a>,<a^{p_1^{\alpha_1}},a^{i_2}b>,<a^{p_2^{\alpha_2}},a^{j_3}b>\},$ \\ $ e_7=\{<a>,<a^{p_1^{\alpha_1}},a^{i_3}b>,<a^{p_2^{\alpha_2}},a^{j_1}b>\}$, where $0 \leq i_1,i_2,i_3 \leq p_1^{\alpha_1}-1$ and $0 \leq j_1,j_2,j_3 \leq  p_2^{\alpha_2}-1$.\\
	Then, the incidence graph $\mathcal{I}(Co_{\mathcal{H}}(D_n))$ contains a subdivision of $K_{3,3}$ as depicted in Figure \ref{planembed}(a). Hence, $Co_{\mathcal{H}}(D_n)$ is not planar.

	\noindent\textbf{Subcase 2.2.}  
	 Suppose $\omega(n)\geq3$, i.e., $n$ has atleast three distinct odd prime divisors of $n$. For the vertices $H_1=<a>, H_2=<a^{p_1},b>, H_3=<a^{p_2},b>$,  $K_1=<a^{p_3},b>, K_2=<a^{p_3},ab>$ and $K_3=<a^{p_3},a^2b>$, where $p_1,p_2,p_3$ are distinct prime divisors of $n$, $p_3$ is an odd prime,  there exist three distinct hyperedges $e_1,e_2,e_3$ of $Co_{\mathcal{H}}(D_n)$ such that $e_1$ contains $H_1,H_2,H_3$ and $K_1$, $e_2$ contains $H_1,H_2,H_3$ and $K_2$, and $e_3$ contains $H_1,H_2,H_3$ and $K_3$. Therefore,  $\mathcal{I}(Co_{\mathcal{H}}(D_n))$ contains $K_{3,3}$. Hence, $Co_{\mathcal{H}}(D_n)$ is non-planar.
		\begin{figure}[h]
		\centering\subfloat[]{	\begin{tikzpicture}[scale=.5, every node/.style={circle, fill=black, minimum size=6pt, inner sep=0pt}]
				
				\node[label=left:${<a>}$] (u1) at (-4,3) {};
				\node[label=left:${<a^{p_1^{\alpha_1}},a^{i_1}b>}$] (u2) at (-4,2) {};
				\node[label=left:${<a^{p_1^{\alpha_1}},a^{i_2}b>}$] (u3) at (-4,1) {};
				\node[label=left:${<a^{p_1^{\alpha_1}},a^{i_3}b>}$] (u4) at (-4,0) {};
				\node[label=left:${<a^{p_2^{\alpha_2}},a^{j_1}b>}$] (u5) at (-4,-1) {};
				\node[label=left:${<a^{p_2^{\alpha_2}},a^{j_2}b>}$] (u6) at (-4,-2) {};
				\node[label=left:${<a^{p_2^{\alpha_2}},a^{j_3}b>}$] (u7) at (-4,-3) {};
				
				\node[label=right:$e_1$] (v1) at (3,3) {};
				\node[label=right:$e_2$] (v2) at (3,2) {};
				\node[label=right:$e_3$] (v3) at (3,1) {};
				\node[label=right:$e_4$] (v4) at (3,0) {};
				\node[label=right:$e_5$] (v5) at (3,-1) {};
				\node[label=right:$e_6$] (v6) at (3,-2) {};
				\node[label=right:$e_7$] (v7) at (3,-3) {};
				\draw (u1) -- (v1);
				\draw (u1) -- (v2);
				\draw (u1) -- (v3);
				\draw (u1) -- (v4);
				\draw (u1) -- (v5);
				\draw (u1) -- (v6);
				\draw (u1) -- (v7);
				\draw (u2) -- (v1);
				\draw (u2) -- (v2);
				\draw (u2) -- (v3);
				\draw (u3) -- (v4);
				\draw (u3) -- (v5);
				\draw (u3) -- (v6);
				\draw (u4) -- (v7);
				\draw (u5) -- (v5);
				\draw (u5) -- (v7);
				\draw (u5) -- (v1);
				\draw (u6) -- (v4);
				\draw (u6) -- (v2);
				\draw (u7) -- (v3);
				\draw (u7) -- (v6);
				
		\end{tikzpicture}}
		\label{G}
		\hspace{.50 cm}
		\subfloat[]{
			\begin{tikzpicture}[scale=.5, every node/.style={circle, fill=black, minimum size=3pt, inner sep=0pt}]
				
				\node [label=left:{$<a^2,ab>$}](u1) at (0,-3) {};
				\node [label=left:$<a>$](u2) at (0,0) {};
				\node [label=left:{$<a^2,b>$}] (u3) at (0,3) {};
				\node [label=left:{}] (u4) at (-2.5,1.5) {};
				\node [label=right:{}] (u5) at (2,-1.5) {};
				\node [label=right:{}] (u6) at (0.3,1.5) {};
				\node [label=right:{}] (u7) at (2.5,1.5) {};
				\node [label=right:{}] (u8) at (3.3,1.5) {};
				\node [label=right:{}] (u9) at (4.3,1.5) {};
				\node [label=right:{}] (u10) at (-4.5,-1.5) {};
				\node [label=right:{}] (u11) at (-0.3,-1.5) {};
				\node [label=right:{}] (u12) at (-3,-1.5) {};
				\node [label=right:{}] (u13) at (-5,-1.5) {};
				\node [label=right:{$e$}](v1) at (7,0) {};
				\node [label=right:{$e_1$}] (v2) at (-1,1.5) {};
				\node [label=left:{$e_1'$}] (v3) at (1,-1.5) {};
				\node [label=right:{$e_2$}] (v4) at (1.5,1.5) {};
				\node [label=above:{$e_3$}] (v5) at (3,1.5) {};
				\node [label=right:{}] (v6) at (3.8,1.5) {};
				\node [label=right:{}] (v7) at (4.8,1.5) {};
				\node [label=left:{$e_2'$}] (v8) at (-1.8,-1.5) {};
				\node [label=above:{$e_3'$}] (v9) at (-4,-1.5) {};
				\node [label=right:{}] (v10) at (-5.5,-1.5) {};
				\node [label=right:{}] (v11) at (-6,-1.5) {};

				\draw (u1) to[out=320, in=270](v1);
				\draw (u2) -- (v1);
				\draw (u3) to[out=20, in=90](v1);
				\draw (u3) -- (v2);
				\draw (u2) -- (v2);
				\draw (u4) -- (v2);
				\draw (u1) -- (v3);
				\draw (u2) -- (v3);
				\draw (u5) -- (v3);
				\draw (u3) -- (v4);
				\draw (u6) -- (v4);
				\draw (u2) -- (v4);
				\draw (u3) -- (v5);
				\draw (u7) -- (v5);
				\draw (u2) -- (v5);
				\draw (u11) -- (v8);
				\draw (u1) -- (v8);
				\draw (u2) -- (v8);
				\draw (u12) -- (v9);
				\draw (u1) -- (v9);
				\draw (u2) -- (v9);
				

		\end{tikzpicture}}
		\caption{}
		\label{planembed}

	\end{figure}

	Conversely, assume that $n$ is a power of a prime. Consider the following cases:\\
	\textbf{Case 1.} If  $n$ is a power of an odd prime, then by Remark \ref{oddprime},  $Co_{\mathcal{H}}(D_n)$ is a 2-uniform star  hypergraph, i.e, a star graph  and hence, $\mathcal{I}(Co_{\mathcal{H}}(D_n))$ can be embedded on a plane. \\ 
		\textbf{Case 2.} If $n$ is a power of 2, then by Theorem \ref{uniformhyp}, $Co_{\mathcal{H}}(D_n)$ is a 3-uniform hypergraph  and hence, $\mathcal{I}(Co_{\mathcal{H}}(D_n))$ can be embedded on the plane as depicted in the Figure \ref{planembed}(b). In the Figure \ref{planembed}, $e, e_1,e_2,e_3 \cdots, e_1', e_2',e_3' \cdots$ are hyperedges of $Co_{\mathcal{H}}(D_n)$ such that $e=\{<a>,<a^2,b>,<a^2,ab>\}$, $e_1,e_2, e_3, \cdots$ are hyperedges containing $<a>, <a^2,b>$ and a \textbf{Type 2} vertex, and  $e_1',e_2', e_3', \cdots$ are hyperedges containing $<a>, <a^2,ab>$ and a \textbf{Type 2} vertex.  Therefore, $Co_{\mathcal{H}}(D_n)$  is planar.
			\end{proof}

\noindent Next, we discuss the possibilities of $Co_{\mathcal{H}}(D_n)$  which can be embedded on the torus and projective plane. The following results about the orientable and non-orientable genus of a hypergraph that are essentially needed to study the embedding of  $Co_{\mathcal{H}}(D_n)$ on these surfaces.
	 
\begin{theorem}\label{incidencegraph}\cite{walsh1975hypermaps}
	For any hypergraph $\mathcal{H}$, $g(\mathcal{H})=g(\mathcal{I}(\mathcal{H}))$.
\end{theorem}
\begin{theorem}\label{incidencegraph1}\cite{walsh1975hypermaps}
	For any hypergraph $\mathcal{H}$, $\tilde{g}(\mathcal{H})=\tilde{g}(\mathcal{I}(\mathcal{H}))$.
\end{theorem}
\begin{lemma}\cite{bojan}\label{genus} The orientable and non-orientable genus of a complete bi-partite graph are given by:
	\begin{enumerate}
		\item $g(K_{m,n})= \left\lceil \frac{(m-2)(n-2)}{4}\right\rceil, m,n \geq 2$
		\item 	$\tilde{g}(K_{m,n}) = 
		\left\lceil \dfrac{(m-2)(n-2)}{2} \right\rceil, \,  m,n\geq 2$
	\end{enumerate}
\end{lemma}

\begin{theorem}\label{g}
	The following statements are equivalent:
	\begin{enumerate}
		\item $n=6$ or $n$ is a power of a prime.
		\item  $Co_{\mathcal{H}}(D_n)$ is toroidal.
		\item $Co_{\mathcal{H}}(D_n)$ is projective.
		
	\end{enumerate}
\end{theorem}
\begin{proof} 1 $\Rightarrow$ 2
	If $n$ is a power of a prime, then by Theorem \ref{planarprimepower}, $Co_{\mathcal{H}}(D_n)$ is planar and hence, it can be embedded on a torus. Therefore, $Co_{\mathcal{H}}(D_n)$ is toroidal. \\Now, 	suppose $n=6$. Then the  hyperedges of $Co_{\mathcal{H}}(D_6)$ are as follows:\\ $e_1=\{<a^2>,<a^3,b>\}, e_2=\{<a^2>,<a^3,ab>\}, $  $e_3=\{<a^2>,<a^3,a^2b>\},$ \\
	$e_4=\{<a^3>,<a^2,ab>,<a^2,b>\}, e_5=\{<b>,<a>,<a^2,ab>\},$ \\$ e_6=\{<a^3b>,<a>,<a^2,b>\},$ $e_7=\{<ab>,<a>,<a^2,b>\},$ 
	$ e_8=\{<a^4b>,<a>,<a^2,ab>\}, $  $e_9=\{<a^2b>,<a>,<a^2,ab>\}, e_{10}=\{<a^5b>,<a>,<a^2,b>\},$ \\ $e_{11}=\{<a>,<a^2,ab>,<a^2,ab>,<a^3,b>\},$  
	$e_{12}=\{<a>,<a^2,ab>,<a^2,ab>,<a^3,ab>\},$ $e_{13}=\{<a>,<a^2,ab>,<a^2,ab>,<a^3,a^2b>\}.$\\
	\noindent The Figure \ref{torusembed} depicts the embedding of  $\mathcal{I}(Co_{\mathcal{H}}(D_6))$ on a torus.  Therefore, $Co_{\mathcal{H}}(D_6)$ is toroidal. 
	\begin{figure}[htbp!]
		\centering
		
		\begin{tikzpicture}[scale=.7, every node/.style={circle, fill=black, minimum size=4pt, inner sep=0pt}]
			\draw (-6,-5) rectangle (7,4); 
			
			\node[label=below:${<a>}$] (u1) at (0,0) {};
			\node[label=left:${<a^2,ab>}$] (u2) at (-6,0) {};
			\node[label=right:${<a^2,b>}$] (u3) at (7,-1) {};
			\node[label=above:${e_5}$] (u4) at (-2,2) {};
			\node[label=left:${<b>}$] (u5) at (-2.5,2.5) {};
			\node[label=above:${e_9}$] (u6) at (-2,1) {};
			\node[label=left:${<a^2b>}$] (u7) at (-.9,0) {};
			\node[label=above:${e_8}$] (u8) at (-2,-2) {};
			\node[label=left:${<a^4b>}$] (u9) at (-1.3,-1) {};
			\node[label=above:${e_7}$] (u10) at (2,2) {};
			\node[label=left:${<ab>}$] (u11) at (3.5,1) {};
			\node[label=above:${e_6}$] (u12) at (2,1) {};
			\node[label=below:${<a^3b>}$] (u13) at (1.7,.5) {};
			\node[label=below:${e_{10}}$] (u14) at (2,-2) {};
			\node[label=right:${<a^5b>}$] (u15) at (1.7,-1) {};
			\node[label=above:${e_{11}}$] (u16) at (-2,4) {};
			\node[label=right:${<a^3,b>}$] (u17) at (-1,3) {};
			\node[label=above left:${e_{1}}$] (u18) at (0,4) {};
			\node[label=left:${<a^2>}$] (u19) at (-2,-4.5) {};
			\node[label=left:${<a^3,ab>}$] (u20) at (-2,-3.5) {};
			\node[label=left:${e_2}$] (u21) at (-1.5,-4) {};
			\node[label=right:${<a^3,a^2b>}$] (u22) at (2,-3.5) {};
			\node[label=below:${e_3}$] (u23) at (0,-4) {};
			\node[label=above:${e_{12}}$] (u24) at (-.5,-3) {};
			\node[label=above:${e_{13}}$] (u25) at (2,4) {};
			\node[label=below:${e_{13}}$] (u26) at (2,-5) {};
			\node[label=right:${<a^2,ab>}$] (u27) at (7,0) {};
			\node[label=below:${e_{1}}$] (u28) at (0,-5) {};
			\node[label=right:${e_{4}}$] (u29) at (6,1) {};
			\node[label=right:$<a^3>$] (u31) at (4.3,2) {};
			\node[label=left:${<a^2,b>}$] (u32) at (-6,-1) {};
			\node[label=below:${e_{11}}$] (u33) at (-2,-5) {};
			\draw (u2) to[out=80, in=180](u4);
			\draw (u4) -- (u5);
			\draw (u4) to[out=300, in=120](u1);
			\draw (u2) to[out=80, in=180](u6);
			\draw (u6) to[out=340, in=180](u1);
			\draw (u6) -- (u7);
			\draw (u2) to[out=360, in=180](u8);
			\draw (u8) to[out=0, in=200](u1);
			\draw (u8) -- (u9);
			\draw (u3) to[out=170, in=0](u10);
			\draw (u10) to[out=180, in=90](u1);
			\draw (u10) -- (u11);
			\draw (u12) -- (u13);
			\draw (u3) to[out=170, in=300](u12);
			\draw (u12) to[out=200, in=10](u1);
			\draw (u1) to[out=330, in=170](u14);
			\draw (u14) to[out=0, in=270](u3);
			\draw (u14) -- (u15);
			\draw (u16) -- (u17);
			\draw (u2) to[out=80, in=180](u16);
			\draw (u16) to[out=270, in=90](u1);
			\draw (u18) -- (u17);
			\draw (u28) to[out=80, in=330](u19);
			\draw (u19) -- (u21);
			\draw (u19) -- (u23);
			\draw (u22) -- (u23);
			\draw (u20) -- (u24);
			\draw (u20) -- (u21);
			\draw (u29) -- (u31);
			\draw (u29) -- (u27);
			\draw (u29) -- (u3);
			
			\draw (u2) to[out=280, in=170](u24);
			\draw (u1) to[out=330, in=30](u24);
			\draw (u32) to[out=300, in=190](u24);
			\draw (u32) to[out=280, in=170](u33);
			\draw  (u26) to[out=90, in=270](u22);
			
			\draw (u27) to[out=100, in=330](u25);
			\draw (u1) to[out=100, in=270](u25);
			\draw (u3) to[out=130, in=300](u25);
			\draw[->,] (0,4) -- (1,4) [right] {};
			\draw[->,] (0,-5) -- (1,-5) [right] {};
			\draw[->,] (7,0) -- (7,1) [right] {};
			\draw[->,] (-6,0) -- (-6,1) [right] {};
		\end{tikzpicture}
		\caption{Embedding of $\mathcal{I}(Co_{\mathcal{H}}(D_6))$ on a torus}
		\label{torusembed}
		
	\end{figure}
	
	\noindent 2 $\Rightarrow$ 1.
 Suppose that $n \neq 6$ and $n$ is not a power of a prime. Then consider the following cases:\\
		\textbf{Case 1.} Suppose $\omega(n)=2$. Consider the following subcases:\\
	\noindent \textbf{Subcase 1.1.} Let $n=2^{\alpha_1}3^{\alpha_2}$, where $\alpha_1, \alpha_2$ are positive integers and atleast one of $\alpha_1,\alpha_2$ is greater than 1. Note that $\mathcal{I}(Co_{\mathcal{H}}(D_6))$ is a subgraph of $\mathcal{I}(Co_{\mathcal{H}}(D_n))$. \\
		\textbf{Subcase 1.1.(a).} If $3^2 \mid n$, then $e'=\{<a>,<a^2,b>,$ $<a^2,ab>,<a^9,b>\}$ and \\$e''=\{<a>,<a^2,b>,$ $<a^2,ab>,<a^9,ab>\}$ are  hyperedges of $Co_{\mathcal{H}}(D_n)$. But there is no way to insert the edges  $\{e',<a>\}, \{e',<a^2,b>\}, \{e',<a^2,ab>\},  \{e',<a^9,b>\}$ $ \{e'',<a>\}, $ \\ $\{e'',<a^2,b>\}, \{e'',<a^2,ab>\}$ and $ \{e'',<a^9,ab>\}$   of $\mathcal{I}(Co_{\mathcal{H}}(D_n))$ in the Figure \ref{torusembed} without crossing. Hence, $\mathcal{I}(Co_{\mathcal{H}}(D_n))$ cannot be embedded on a torus and consequently, $Co_{\mathcal{H}}(D_n)$ is not toroidal.\\
	\textbf{Subcase 1.1.(b).} If $2^2\mid n$, then $e'''=\{<a>,<a^4,b>,<a^3,b>,<a^2,ab>\}$ is a hyperedge of $Co_{\mathcal{H}}(D_n)$. But there is no way to insert the edges $\{e''',<a>\}, \{e''',<a^3,b>\}, \{e''',<a^2,ab>\}, $ and $ \{e''',<a^4,b>\}$ of $\mathcal{I}(Co_{\mathcal{H}}(D_n))$ in the Figure \ref{torusembed} without crossing. Hence, $\mathcal{I}(Co_{\mathcal{H}}(D_n))$ cannot be embedded on a torus. Thus, $Co_{\mathcal{H}}(D_n)$ is not toroidal.\\
	\noindent \textbf{Subcase 1.2.}  Let $n=2^{\alpha_1}5^{\alpha_2}$, where $\alpha_1, \alpha_2$ are positive integers.	Consider a subhypergraph $G$ of $Co_{\mathcal{H}}(D_n)$ whose hyperedges are as follows:\\ $e_1=\{<a^2>,<a^5,b>\}, \, e_2=\{<a^2>,<a^5,ab>\}, $  $e_3=\{<a^2>,<a^5,a^2b>\},$ \\ $ e_4=\{<a^3>,<a^2,ab>,<a^2,b>\},\, e_5=\{<b>,<a>,<a^2,ab>\},$ \\ $e_6=\{<a^3b>,<a>,<a^2,b>\},  e_7=\{<ab>,<a>,<a^2,b>\}, \, e_8=\{<a^4b>,<a>,<a^2,ab>\},$ \\$ e_9=\{<a^2b>,<a>,<a^2,ab>\}, \, e_{10}=\{<a^5b>,<a>,<a^2,b>\},$ \\ $ e_{11}=\{<a>,<a^2,ab>,<a^2,ab>,<a^5,b>\}, \, e_{12}=\{<a>,<a^2,ab>,<a^2,ab>,<a^5,ab>\},$ $ e_{13}=\{<a>,<a^2,ab>,<a^2,ab>,<a^5,a^2b>\}, \,  e_{14}=\{<a^2>,<a^5,a^3b>\}, $ \\ $ e_{15}=\{<a>,<a^5,a^3b>,<a^2,b>a^2,ab>\}, \, e_{16}=\{<a>,<a^5,a^4b>,<a^2,b>,<a^2,ab>\}$. 
	 \begin{figure}[htbp!]
		\centering
		
		\begin{tikzpicture}[scale=.8, every node/.style={circle, fill=black, minimum size=4pt, inner sep=0pt}]
			\draw (-6,-5) rectangle (7,4); 
			
			\node[label=below:${ <a>}$] (u1) at (0,0) {};
			\node[label=left:${<a^2,ab>}$] (u2) at (-6,0) {};
			\node[label=right:${<a^2,b>}$] (u3) at (7,-1) {};
			\node[label=above:${e_5}$] (u4) at (-2,2) {};
			\node[label=left:${<b>}$] (u5) at (-2.5,2.5) {};
			\node[label=above:${e_9}$] (u6) at (-2,1) {};
			\node[label=left:${<a^2b>}$] (u7) at (-.9,0) {};
			\node[label=above:${e_8}$] (u8) at (-2,-2) {};
			\node[label=left:${<a^4b>}$] (u9) at (-1.3,-1) {};
			\node[label=above:${e_7}$] (u10) at (2,2) {};
			\node[label=left:${<ab>}$] (u11) at (3.5,1) {};
			\node[label=above:${e_6}$] (u12) at (2,1) {};
			\node[label=below:${<a^3b>}$] (u13) at (1.7,.5) {};
			\node[label=below:${e_{10}}$] (u14) at (2,-2) {};
			\node[label=right:${<a^5b>}$] (u15) at (1.7,-1) {};
			\node[label=above:${e_{11}}$] (u16) at (-2,4) {};
			\node[label=right:${<a^5,b>}$] (u17) at (-1,3) {};
			\node[label=above left:${e_{1}}$] (u18) at (0,4) {};
			\node[label=left:${<a^2>}$] (u19) at (-2,-4.5) {};
			\node[label=left:${<a^5,ab>}$] (u20) at (-2,-3.5) {};
			\node[label=left:${e_2}$] (u21) at (-1.5,-4) {};
			\node[label=right:\scriptsize{${<a^5,a^2b>}$}] (u22) at (2,-3.5) {};
			\node[label=below:${e_3}$] (u23) at (0,-4) {};
			\node[label=above:${e_{12}}$] (u24) at (-.5,-3) {};
			\node[label=above:${e_{13}}$] (u25) at (2,4) {};
			\node[label=below:${e_{13}}$] (u26) at (2,-5) {};
			\node[label=right:${<a^2,ab>}$] (u27) at (7,0) {};
			\node[label=below:${e_{1}}$] (u28) at (0,-5) {};
			\node[label=right:${e_{4}}$] (u29) at (6,1) {};
			\node[label=right:$<a^5>$] (u31) at (4.3,2) {};
			\node[label=left:${<a^2,b>}$] (u32) at (-6,-1) {};
			\node[label=below:${e_{11}}$] (u33) at (-2,-5) {};
			\node[label=above right:${e_{14}}$] (u34) at (0,-3) {};
			\node[label=right:\scriptsize{${<a^5,a^3b>}$}] (u35) at (4,-3.5) {};
			\node[label=below:${e_{15}}$] (u36) at (5,-5) {};
			\node[label=above:${e_{15}}$] (u37) at (5,4) {};
			\draw (u2) to[out=80, in=180](u4);
			\draw (u4) -- (u5);
			\draw (u4) to[out=300, in=120](u1);
			\draw (u2) to[out=80, in=180](u6);
			\draw (u6) to[out=340, in=180](u1);
			\draw (u6) -- (u7);
			\draw (u2) to[out=360, in=180](u8);
			\draw (u8) to[out=0, in=200](u1);
			\draw (u8) -- (u9);
			\draw (u3) to[out=170, in=0](u10);
			\draw (u10) to[out=180, in=90](u1);
			\draw (u10) -- (u11);
			\draw (u12) -- (u13);
			\draw (u3) to[out=170, in=300](u12);
			\draw (u12) to[out=200, in=10](u1);
			\draw (u1) to[out=330, in=170](u14);
			\draw (u14) to[out=0, in=270](u3);
			\draw (u14) -- (u15);
			\draw (u16) -- (u17);
			\draw (u2) to[out=80, in=180](u16);
			\draw (u16) to[out=270, in=90](u1);
			\draw (u18) -- (u17);
			\draw (u28) to[out=80, in=330](u19);
			\draw (u19) -- (u21);
			\draw (u19) -- (u23);
			\draw (u22) -- (u23);
			\draw (u20) -- (u24);
			\draw (u20) -- (u21);
			\draw (u29) -- (u31);
			\draw (u29) -- (u27);
			\draw (u29) -- (u3);
			\draw (u36) -- (u35);
			
			\draw (u2) to[out=280, in=170](u24);
			\draw (u1) to[out=240, in=80](u24);
			\draw (u32) to[out=300, in=190](u24);
			\draw (u32) to[out=280, in=170](u33);
			\draw  (u26) to[out=90, in=270](u22);
			\draw  (u26) to[out=90, in=270](u22);
			\draw (u1) .. controls (0,-3) and (6,-2) .. (u36);
			\draw (u27) to[out=100, in=330](u25);
			\draw (u1) to[out=100, in=270](u25);
			\draw (u3) to[out=130, in=300](u25);
			\draw (u19) to[out=30, in=270](u34);
			\draw (u34) to[out=330, in=90](u35);
			\draw (u3) to[out=270, in=30](u36);
			\draw (u27) to[out=100, in=320](u37);
			\draw[->,] (0,4) -- (1,4) [right] {};
			\draw[->,] (0,-5) -- (1,-5) [right] {};
			\draw[->,] (7,0) -- (7,1) [right] {};
			\draw[->,] (-6,0) -- (-6,1) [right] {};
		\end{tikzpicture}
		\caption{}
		\label{torusembed3}
		
	\end{figure}
	
	\noindent The Figure \ref{torusembed3} depicts the embedding of a subgraph of   $\mathcal{I}(G)$ on a torus.  Inserting the edges $\{e_{16},<a>\},$ $ \{e_{16},<a^5,a^4b>\}, \{e_{16},<a^2,b>\}, \{e_{16},<a^2,ab>\}$ of $\mathcal{I}(G)$ in the Figure \ref{torusembed3}   without crossing is not possible. Therefore, $Co_{\mathcal{H}}(D_n)$ is not toroidal.\\
	\textbf{Subcase 1.3.} Let  $n=3^{\alpha_1}5^{\alpha_2}$,  where $\alpha_1, \alpha_2$ are positive integers. Consider the subhypergraph $G'$ of $Co_{\mathcal{H}}(D_n)$, where the hyperedge set of $G'$ are as follows: $e_1=\{<a^3>,<a^5,b>\},$ \\ $ e_2=\{<a^3>,<a^5,ab>\}, e_3=\{<a^3>,<a^5,a^2b>\},$ $ e_4=\{<a^3>,<a^5,a^3b>\},$ \\ $ e_5=\{<a^3>,<a^5,a^4b>\}, e_6=\{<a^5>,<a^3,b>\},  e_7=\{<a^5>,<a^3,ab>\},$ \\ $ e_8=\{<a^5>,<a^3,a^2b>\}, e_9=\{<a>,<b>\}, e_{10}=\{<a>,<ab>\}, e_{11}=\{<a>,<a^2b>\},$ $ e_{12}=\{<a>,<a^3b>\}, e_{13}=\{<a>,<a^4b>\},$ $ e_{14}=\{<a>,<a^5b>\}, e_{15}=\{<a>,<a^6b>\},$ $ e_{16}=\{<a>,<a^7b>\},e_{17}=\{<a>,<a^8b>\}, e_{18}=\{<a>,<a^9b>\}, e_{19}=\{<a>,<a^{10}b>\},$ $ e_{20}=\{<a>,<a^{11}b>\}, e_{21}=\{<a>,<a^{12}b>\},  e_{22}=\{<a>,<a^{13}b>\}, e_{23}=\{<a>,<a^{14}b>\},$ $ e_{24}=\{<a>,<a^3,b>,<a^5,b>\}, e_{25}=\{<a>,<a^3,b>,<a^5,ab>\},$ \\ $ e_{26}=\{<a>,<a^3,b>,<a^5,a^2b>\}, e_{27}=\{<a>,<a^3,b>,<a^5,a^3b>\},$\\ $e_{28}=\{<a>,<a^3,b>,<a^5,a^4b>\}$. The vertex set of $G'$ is all those vertices in hyperedges of $e_1,e_2, e_3, \cdots$ and $e_{28}$. The figure \ref{torus2} depicts the embedding of  $\mathcal{I}(G') - \{e_{28},<a>\}$ on a torus and we cannot insert the edge \{$e_{28},<a>$\} without crossing. Hence, $Co_{\mathcal{H}}(D_n)$ is not toroidal.
	\begin{figure}[htbp!]
		\centering
		
		\begin{tikzpicture}[scale=.75, every node/.style={circle, fill=black, minimum size=4pt, inner sep=0pt}]
			\draw (-9,-5) rectangle (7,4); 
			
			\node[label=below:${<a>}$] (u1) at (0,-5) {};
			\node[label=right:${<a^3>}$] (u2) at (-4,0) {};
			\node[label=left:${<a^5,b>}$] (u3) at (-6.3,2) {};
			\node[label=left:${<a^5,ab>}$] (u4) at (-6.3,1) {};
			\node[label=left:${<a^5,a^2b>}$] (u5) at (-6.3,0) {};
			\node[label=left:${<a^5,a^3b>}$] (u6) at (-6.3,-1) {};
			\node[label=left:${<a^5,a^4b>}$] (u7) at (-6.3,-2) {};
			\node[label=right:${<a^5>}$] (u8) at (4.5,0) {};
			\node[label=left:${<a^3,b>}$] (u9) at (3.5,1) {};
			\node[label=left:${<a^3,ab>}$] (u10) at (3.5,0) {};
			\node[label=left:${<a^3,a^2b>}$] (u11) at (3.5,-1) {};
			\node[label=right:${e_1}$] (u12) at (-6,3) {};
			\node[label=left:${e_2}$] (u13) at (-5.6,2) {};
			\node[label=above:${e_3}$] (u14) at (-5.4,1) {};
			\node[label=above:${e_4}$] (u15) at (-5.4,0) {};
			\node[label=above:${e_5}$] (u16) at (-5.4,-1) {};
			\node[label=above:${e_6}$] (u17) at (4.5,2) {};
			\node[label=above:${e_7}$] (u18) at (4.2,.8) {};
			\node[label=below:${e_8}$] (u19) at (4.2,-.8) {};
			\node[label=left:${e_9}$] (u20) at (-5.8,-4) {};
			\node[label={[rotate=90]right:\tiny{$<b>$}}] (u21) at (-5.8,-3.5) {};
			\node[label=left:${e_{10}}$] (u22) at (-5.0,-4) {};
			\node[label={[rotate=90]right:\tiny{$<ab>$}}] (u23) at (-5,-3.5) {};
			\node[label=left:${e_{11}}$] (u24) at (-4.2,-4) {};
			\node[label={[rotate=90]right:\tiny{$<a^2b>$}}] (u25) at (-4.2,-3.5) {};
			\node[label=left:${e_{12}}$] (u26) at (-3.4,-4) {};
			\node[label=left:${e_{13}}$] (u27) at (-2.6,-4) {};
			\node[label=left:${e_{14}}$] (u28) at (-1.8,-4) {};
			\node[label=left:${e_{15}}$] (u29) at (-1,-4) {};
			\node[label=left:${e_{16}}$] (u30) at (-0.2,-4) {};
			\node[label=left:${e_{17}}$] (u31) at (0.6,-4) {};
			\node[label=left:${e_{18}}$] (u32) at (1.4,-4) {};
			\node[label=left:${e_{19}}$] (u33) at (2.2,-4) {};
			\node[label=left:${e_{20}}$] (u34) at (3,-4) {};
			\node[label=left:${e_{21}}$] (u35) at (3.8,-4) {};
			\node[label=left:${e_{22}}$] (u36) at (4.6,-4) {};
			\node[label=left:${e_{23}}$] (u37) at (5.4,-4) {};
			\node[label={[rotate=90]right:\tiny{$<a^3b>$}}] (u38) at (-3.4,-3.5) {};
			\node[label={[rotate=90]right:\tiny{$<a^4b>$}}] (u39) at (-2.6,-3.5) {};
			\node[label={[rotate=90]right:\tiny{$<a^5b>$}}] (u40) at (-1.8,-3.5) {};
			\node[label={[rotate=90]right:\tiny{$<a^6b>$}}] (u41) at (-1,-3.5) {};
			\node[label={[rotate=90]right:\tiny{$<a^7b>$}}] (u42) at (-.2,-3.5) {};
			\node[label={[rotate=90]right:\tiny{$<a^8b>$}}] (u43) at (0.6,-3.5) {};
			\node[label={[rotate=90]right:\tiny{$<a^9b>$}}] (u44) at (1.4,-3.5) {};
			\node[label={[rotate=90]right:\tiny{$<a^{10}b>$}}] (u45) at (2.2,-3.5) {};
			\node[label={[rotate=90]right:\tiny{$<a^{11}b>$}}] (u46) at (3,-3.5) {};
			\node[label={[rotate=90]right:\tiny{$<a^{12}b>$}}] (u47) at (3.8,-3.5) {};
			\node[label={[rotate=90]right:\tiny{$<a^{13}b>$}}] (u48) at (4.6,-3.5) {};
			\node[label={[rotate=90]right:\tiny{$<a^{14}b>$}}] (u49) at (5.4,-3.5) {};
			\node[label=below left:${e_{24}}$] (u50) at (-1,3) {};
			\node[label=below:${<a>}$] (u51) at (0,4) {};
			\node[label=left:${e_{25}}$] (u52) at (-9,0) {};
			\node[label=right:${e_{25}}$] (u53) at (7,0) {};
			\node[label=left:${e_{26}}$] (u54) at (-9,-1) {};
			\node[label=left:${e_{26}}$] (u55) at (7,-1) {};
			\node[label=left:${e_{27}}$] (u56) at (-1,1) {};
			\node[label=left:${e_{28}}$] (u57) at (-1,2) {};

			\draw (u20) to[out=310, in=170](u1);
			\draw (u12) -- (u2);
			\draw (u12) -- (u3);
			\draw (u13) -- (u2);
			\draw (u13) -- (u4);
			\draw (u14) -- (u2);
			\draw (u14) -- (u5);
			\draw (u15) -- (u2);
			\draw (u15) -- (u6);
			\draw (u16) -- (u2);
			\draw (u16) -- (u7);
			\draw (u17) -- (u9);
			\draw (u17) -- (u8);
			\draw (u18) -- (u10);
			\draw (u18) -- (u8);
			\draw (u19) -- (u11);
			\draw (u19) -- (u8);
			\draw (u21) -- (u20);
			\draw (u22) -- (u23);
			\draw (u24) -- (u25);
			\draw (u26) -- (u38);
			\draw (u27) -- (u39);
			\draw (u28) -- (u40);
			\draw (u29) -- (u41);
			\draw (u30) -- (u42);
			\draw (u31) -- (u43);
			\draw (u32) -- (u44);
			\draw (u33) -- (u45);
			\draw (u34) -- (u46);
			\draw (u35) -- (u47);
			\draw (u36) -- (u48);
			\draw (u37) -- (u49);
			\draw (u51) -- (u50);
			\draw (u9) -- (u50);
			
			\draw (u22) to[out=310, in=170](u1);
			\draw (u24) to[out=310, in=170](u1);
			\draw (u26) to[out=310, in=170](u1);
			\draw (u27) to[out=310, in=170](u1);
			\draw (u28) to[out=310, in=170](u1);
			\draw (u29) to[out=310, in=170](u1);
			\draw (u30) to[out=310, in=170](u1);
			\draw (u31) to[out=310, in=10](u1);
			\draw (u32) to[out=310, in=10](u1);
			\draw (u33) to[out=310, in=10](u1);
			\draw (u34) to[out=310, in=10](u1);
			\draw (u35) to[out=310, in=10](u1);
			\draw (u36) to[out=310, in=10](u1);
			\draw (u37) to[out=310, in=10](u1);
			\draw (u3) to[out=120, in=170](u50);
			\draw (u52) to[out=80, in=185](u51);
			\draw (u52) to[out=80, in=260](u4);
			\draw (u54) to[out=70, in=200](u5);
			\draw (u54) .. controls (-8.8,-5) and (1.2,0) .. (u1);
			\draw (u55) .. controls (4,4) .. (u9);
			\draw (u53) .. controls  (4,4.5) .. (u9);
			\draw (u56) .. controls (0,2) .. (u9);
			\draw (u6) .. controls (-7,-3) .. (u56);
			\draw (u56) .. controls (0.5,-4) .. (u1);
			\draw (u57) .. controls (-3,0) and (-.3,2) .. (u7);
			\draw (u57).. controls (0,2.5) .. (u9);
			\draw[->,] (0,4) -- (1,4) [right] {};
			\draw[->,] (0,-5) -- (1,-5) [right] {};
			\draw[->,] (7,0) -- (7,1) [right] {};
			\draw[->,] (-9,0) -- (-9,1) [right] {};
		\end{tikzpicture}
		\caption{Embedding of $\mathcal{I}(G') - \{e_{28},<a>\}$ on a torus}
		\label{torus2}
		
	\end{figure}

	\noindent \textbf{Subcase 1.4.} Let  $n=p_1^{\alpha_1} p_2^{\alpha_2}$, where $p_1,p_2$ are odd primes, $p_2 \geq 7$ and  $\alpha_1, \alpha_2$ are positive integers. Then, for the vertices $H_1=<a>, H_2=<a^{p_1},b>,   H_3=<a^{p_1},ab>, K_1=<a^{p_2},b>, K_2=<a^{p_2},ab>,$ $K_3=<a^{p_2},a^2b>,$ $ K_4=<a^{p_2},a^3b>, K_5=<a^{p_2},a^4b>$,  $K_6=<a^{p_2},a^5b>$ and $K_7=<a^{p_2},a^6b>$,  there exist seven distinct hyperedges $e_1,e_2,e_3, e_4,e_5,e_6$ and $e_7$ such that $e_1$ contains $H_1,H_2,H_3$ and $K_1$, $e_2$ contains $H_1,H_2,H_3$ and $K_2$, $e_3$ contains $H_1,H_2,H_3$ and $K_3$, $e_4$ contains $H_1,H_2,H_3$ and $K_4$, $e_5$ contains $H_1,H_2,H_3$ and $K_5$, $e_6$ contains $H_1,H_2,H_3$ and $K_6$ and $e_7$ contains $H_1,H_2,H_3$ and $K_7$. Thus, the incidence graph $\mathcal{I}(Co_{\mathcal{H}}(D_n))$ contains $K_{3,7}$ as a subgraph. By Lemma \ref{genus}, $g(K_{3,7}) = 2$. Consequently, by Theorem \ref{incidencegraph},  $g(Co_{\mathcal{H}}(D_n)) \geq 2.$ Therefore, $Co_{\mathcal{H}}(D_n)$ is not toroidal.\\
	\noindent \textbf{Case 2.} Suppose $\omega(n)\geq 3$. Consider the following subcases:\\
		\noindent \textbf{Subcase 2.1.} Suppose $n=2^{\alpha_1} 3^{\alpha_2}5^{\alpha_3}$, where $\alpha_1, \alpha_2$ are positive integers. Consider the hyperedges $e_,e_2, e_3,e_4,e_5,e_6$  of $Co_{\mathcal{H}}(D_n)$ as follows: $e_1 \supseteq \{<a>,<a^2,b>,<a^3,b>,<a^5,b>\},$ \\$  e_2 \supseteq \{<a>,<a^2,b>,<a^3,b>,<a^5,ab>\},$ $ e_3 \supseteq \{<a>,<a^2,b>,<a^3,b>,<a^5,a^3b>\},$ $ e_4 \supseteq \{<a>,<a^2,ab>,<a^3,ab>,<a^5,b>\},$ $ e_5 \supseteq \{<a>,<a^2,ab>,<a^3,ab>,<a^5,ab>\}$, $ e_6 \supseteq \{<a>,<a^2,ab>,<a^3,ab>,<a^5,a^3b>\}$. Then, $\mathcal{I}(Co_{\mathcal{H}}(D_n))$ contains a subgraph  as shown in Figure \ref{torusobstruction} which is a torus obstruction. Hence, $Co_{\mathcal{H}}(D_n)$ cannot be embedded on a torus. Therefore, $Co_{\mathcal{H}}(D_n)$ is not toroidal.
			\begin{figure}[htbp!]
			\centering
			
			\begin{tikzpicture}[scale=1.5, rotate=90,  every node/.style={circle, fill=black, minimum size=6pt, inner sep=0pt}]
				
				\node[label=left:{$<a^2,b>$}] (u1) at (0,3) {};
				\node[label=right:{$<a>$}] (u2) at (-0.5,2) {};
				\node[label=right:{$<a^3,b>$}] (u3) at (0,1) {};
				\node[label=left:{$<a^2,ab>$}] (u4) at (-.5,3) {};
				\node[label=right:{$<a^3,ab>$}] (u5) at (-.5,1) {};
				\node[label=right:$e_1$] (v1) at (1,3) {};
				\node[label=right:$e_2$] (v2) at (1,2) {};
				\node[label=right:$e_3$] (v3) at (1,1) {};
				\node[label=left:$e_4$] (v4) at (-1.5,3) {};
				\node[label=left:$e_5$] (v5) at (-1.5,2) {};
				\node[label=left:$e_6$] (v6) at (-1.5,1) {};
				\draw (u1) -- (v1);
				\draw (u1) -- (v2);
				\draw (u1) -- (v3);
				\draw (u2) -- (v1);
				\draw (u2) -- (v2);
				\draw (u2) -- (v3);
				\draw (u3) -- (v3);
				\draw (u3) -- (v2);
				\draw (u3) -- (v1);
				\draw (u4) -- (v4);
				\draw (u4) -- (v5);
				\draw (u4) -- (v6);
				\draw (u5) -- (v4);
				\draw (u5) -- (v5);
				\draw (u5) -- (v6);
				\draw (u2) -- (v4);
				\draw (u2) -- (v5);
				\draw (u2) -- (v6);
			\end{tikzpicture}
			\caption{Subgraph of $\mathcal{I}(Co_\mathcal{H}(D_n))$}
			\label{torusobstruction}
		\end{figure}
				
	\noindent \textbf{Subcase 2.2.} Suppose that $n$ has a prime divisor greater than or equal to 7. 
 Let  $p_1,p_2$ and $p_3$ be distinct prime divisors of $n$. Without loss of generality, assume that $p_3 \geq 7$. Then, for the vertices $H_1=<a>, H_2=<a^{p_1},b>, H_3=<a^{p_2},b>$, $K_1=<a^{p_3},b>, K_2=<a^{p_3},ab>, K_3=<a^{p_3},a^2b>,$ $K_4= <a^{p_3},a^3b>,$ $ K_5=<a^{p_3},a^4b>$,  $K_6=<a^{p_3},a^5b>$ and $K_7=<a^{p_3},a^6b>$,  there exist seven distinct hyperedges $e_1,e_2,e_3, e_4,e_5,e_6$ and $e_7$ such that $e_1$ contains $H_1,H_2,H_3$ and $K_1$, $e_2$ contains $H_1,H_2,H_3$ and $K_2$, $e_3$ contains $H_1,H_2,H_3$ and $K_3$, $e_4$ contains $H_1,H_2,H_3$ and $K_4$, $e_5$ contains $H_1,H_2,H_3$ and $K_5$, $e_6$ contains $H_1,H_2,H_3$ and $K_6$ and $e_7$ contains $H_1,H_2,H_3$ and $K_7$. Thus, the incidence graph $\mathcal{I}(Co_{\mathcal{H}}(D_n))$ contains $K_{3,7}$ as a subgraph. By Lemma \ref{genus}, $g(K_{3,7}) = 2$. Consequently, by Theorem \ref{incidencegraph},  $g(Co_{\mathcal{H}}(D_n)) \geq 2.$ Therefore, $Co_{\mathcal{H}}(D_n)$ is not toroidal.

\noindent By similar arguments,  $1 \iff 3$ except  the Subcase 2.1., i.e., for $n=2^{\alpha_1} 3^{\alpha_2}5^{\alpha_3}$, where $\alpha_1, \alpha_2$ are positive integers.  In this case, for the vertices $H_1=<a>,  H_2=<a^{2},b>, H_3=<a^{2},ab>,$\\$ K_1=<a^{5},b>, K_2=<a^{5},ab>, K_3=<a^{5},a^2b>,$ $ K_4=<a^{5},a^3b>, K_5=<a^{5},a^4b>$,  there exist five distinct hyperedges $e_1,e_2,e_3, e_4$ and $e_5$ such that $e_1$ contains $H_1,H_2,H_3$ and $K_1$, $e_2$ contains $H_1,H_2,H_3$ and $K_2$, $e_3$ contains $H_1,H_2,H_3$ and $K_3$, $e_4$ contains $H_1,H_2,H_3$ and $K_4$, $e_5$ contains $H_1,H_2,H_3$ and $K_5$.  Thus, the incidence graph $\mathcal{I}(Co_{\mathcal{H}}(D_n))$  contains $K_{3,5}$ as a subgraph. By Lemma \ref{genus}, $\tilde{g}(K_{3,5}) = 2$. Consequently, by Theorem \ref{incidencegraph1},  $\tilde{g}(Co_{\mathcal{H}}(D_n)) \geq 2$ and so, $Co_{\mathcal{H}}(D_n)$ is not projective.
\end{proof}

\section*{Concluding Remark}
The study of hypergraphs associated with the subgroups of a group is an emerging area of research. In this paper, we studied the co-maximal subgroups of dihedral groups from the hypergraph perspective. The study of co-maximal hypergraphs of a group is nothing but the study of maximal cliques of the deleted co-maximal subgroup graph of a group. So, we have identified some structural properties of deleted co-maximal hypergraphs which are invariant and which are not.

\section*{Acknowledgment}
This research was financially supported by the University Grants Commission (UGC), Government of India, through the award of the Junior Research Fellowship (JRF) under the National Eligibility Test (NET) (Ref. No. 241610038497) awarded to the second author.

\end{document}